\documentclass{article}

\usepackage{amsfonts,amsmath}
\usepackage{graphicx}        % standard LaTeX graphics tool
\newtheorem{prop}{Proposition}

\newtheorem{defn}{Definition}
\newtheorem{rem}{Remark}

\newcommand{\RR}{{\mathbb R}}

\newcommand{\norm}[1]{\lVert#1\rVert}

\newcommand{\Rset}{{\mathbb R}}

\newcommand{\XX}{{\mathcal X }}

\newcommand{\UU}{{\mathcal  U}}
\newcommand{\LL}{{\mathcal  L}}

\newcommand{\dv}[2]{{\frac{\partial #1}{\partial #2}}}

\newcommand{\dotex}{{\frac{d}{dt}}}

\begin{document}

\title{Symmetries in observer design: review of some recent results and applications to EKF-based SLAM}
\author{Silv\`{e}re Bonnabel\footnote{Silv\`{e}re Bonnabel, Centre de Robotique,
Math\'{e}matiques et Syst\`{e}mes, Mines ParisTech, 60 bd Saint-Michel, 75272 Paris cedex 06. {silvere.bonnabel@mines-paristech.fr}}
% Use \authorrunning{Short Title} for an abbreviated version of
% your contribution title if the original one is too long
%\institute{Silv\`{e}re Bonnabel, Centre de Robotique,
%Math\'{e}matiques et Syst\`{e}mes, Mines ParisTech, 60 bd Saint-Michel, 75272 Paris %cedex 06. \email{silvere.bonnabel@mines-paristech.fr}
}

\maketitle

\abstract{In this paper, we first review the theory of
symmetry-preserving observers and we mention some recent results.
Then, we apply the theory to Extended Kalman Filter-based
Simultaneous Localization and Mapping (EKF SLAM). It allows to
derive  a new  (symmetry-preserving) Extended Kalman Filter for the
non-linear SLAM problem that possesses  convergence properties. We
also prove a special choice of the gains ensures global
exponential convergence. }

\section{Introduction}

Symmetries and Lie groups have  been widely used for feedback control in robotics, see e.g. \cite{bullo-murray-auto99,morin-ieee03}.  More generally control of systems possessing symmetries has also been studied for quite a long time, see e.g. \cite{grizzle-marcus-ieee85,martin-et-al-cocv03}. The use of symmetries and Lie groups for observer design is more recent \cite{aghannan-rouchon-cdc02,arxiv-07}. The main properties of those observers are based on the reduction of the  estimation error complexity. When the symmetry group coincides with the state space (observers on Lie groups), the error equation can be particularly simple \cite{arxiv-08}. This property has been used to derive non-linear observers with (almost) global convergence properties for several  localisation problems \cite{mahony-et-al-IEEE,arxiv-08,vasconcelos}.  Recently \cite{bonnabel-et-al:ifac11} established a link between observer design and control of systems on Lie group by proving a non-linear separation principle on Lie groups. 

This paper proposes to recap the main elements of the theory along with some recent results, and to apply it to  the domain of Extended Kalman Filter-based Simultaneous Localization and Mapping (EKF SLAM). It is organized as follows: Section 2 is a brief recap on linear observers. In Section 3 we recap  the theory of symmetry-preserving observers \cite{arxiv-07} and mention some recent results \cite{bonnabel-et-al:ifac11,bonnabel-et-al:CDC11}. In Section 4 we apply it in a straightforward way to EKF SLAM. In Section 5 some results for the special case of observers for invariant systems on Lie groups \cite{arxiv-08} are recalled.  In Section 6, it is proved that those results can be  (surprisingly) applied to EKF SLAM. We  derive a simple globally convergent observer for the non-linear problem. We also propose a modified EKF  such that the covariance matrix and the gain matrix behave as if the system was linear and time-invariant. Such non-linear convergence guarantees for EKF SLAM are new to the author's knowledge.
The author would like to mention and to thank his regular co-authors on the subject of symmetry-preserving observers : Philippe Martin, Pierre Rouchon, and Erwan Sala\"{u}n.

\section{Luenberger observers, extended Kalman filter}
\subsection{Observers for linear systems}
\label{sec:linear_observers}
Observers are meant to compute an estimation of the state of a dynamical system from several sensor measurements. Let $x\in \RR^n$ denote the state of the system, $u\in\RR^m$ be the inputs (a set of $m$ known scalar variables such as controls, constant parameters, etc.). We assume the sensors provide measurements $y\in\RR^p$ that can be expressed as a function of the state and the inputs. When the underlying dynamical model is a linear differential equation, and the output is a linear function as well, the system can be written
\begin{equation}\label{eq:lin_sys}
  \dotex x=A x+Bu, \qquad
 y=Cx+Du.\end{equation}A Luenberger observer (or Kalman filter) writes \begin{align}\label{lin:obs:eq}
  \dotex \hat x=A \hat x+Bu-{L\cdot(C\hat x+Du-y)},\end{align}
  where $\hat x$ is the estimated state, and L is a {gain} matrix that can be freely chosen. We see that the observer consists in  a copy of the system dynamics $A \hat x+Bu$, plus a correction term ${L(C\hat x+Du-y)}$ that ``corrects" the trusted dynamics in function of the discrepancy between the estimated output ${\hat y=C\hat x+Du}$ and the measured output $y$.

One important issue is the choice (or ``tuning") of the gain matrix $L$. The Luenberger observer is based on a choice of a fixed matrix $L$. In the Kalman filter two positive definite matrices $M$ and $N$  denote the covariance matrices of the state noise and measurement noise, and $L$ relies on a Ricatti equation : $
 L =PC^TN,$ where $
\dotex P =AP+PA^T+M^{-1}-PC^TNCP.$
As $M$ and $N$ must be defined by the user,   they can be viewed as tuning matrices.

 In both cases the observer has the form \eqref{lin:obs:eq} with $L$ constant or not. Let $\tilde x=\hat x-x$ be the estimation error, and let us compute the differential equation  satisfied by the error. We have
\begin{align}\label{err:eq3}
\dotex \tilde x= (A+LC)\tilde x.
\end{align}
As the goal of the observer is to find an estimate of $x$, we want $\tilde x$ to go to zero. When the system is observable, one can always find $L$ such that $\tilde x$ asymptotically exponentially goes to zero, and the negative real part of the eigenvalues of  $A+LC$ can be freely assigned. We see that the theory is particularly simple as the error equation~\eqref{err:eq3} is \emph{autonomous}, \emph{i.e.} it does not depend on the trajectory followed by the system. In particular, the input term $u$ has vanished in~\eqref{err:eq3}. The well-known separation principle stems from this fact.

\subsection{Some popular extensions to nonlinear systems}
\label{sec:nonlinear_observers}
Consider a general nonlinear system
\begin{align}\label{nonlin:eq1}
 \dotex x =f(x,u), \qquad y =h(x,u),
 \end{align}
where
  $x\in\mathcal X\subset\mathbb{R}^n$ is the state, $u\in\mathcal U\subset\mathbb{R}^m$
  the input, and $y\in\mathcal Y\subset\mathbb{R}^p$ the output. Mimicking the linear case, a class of popular nonlinear observers writes
  \begin{align}\label{asymp:obs:eq}\dotex{\hat x}=f(\hat x,u)-L(\hat x,y,t)\cdot\bigl(h(\hat x,u)-y(t)\bigr),\end{align}where the gain matrix can depend on the variables $\hat x,y,t$. The error equation can still be computed, but as the system is nonlinear, it does not necessarily lead to an appropriate  gain matrix $L$. Indeed we have $
  \dotex{\tilde x}=f(\hat x,u(t))-f( x,u(t))-L(\hat x,y(t),t)\cdot\bigl(h(\hat x,u(t))-y(t)\bigr).
  $ The error equation is no longer autonomous, and the problem of finding $L$ such that $\tilde x$ goes asymptotically to zero can not be solved in the general  case.

  The most popular  observer for  nonlinear systems
is the Extended Kalman Filter (EKF). The principle is to linearize the system around the estimated trajectory, build a Kalman filter for the linear model, and implement it on the nonlinear system.  The EKF has the form  \eqref{asymp:obs:eq}, where the gain matrix is computed the following way:  \begin{align}\label{kalman:gain1}
 A &=\dv{f}{x}(\hat x,u) &\quad  C &=\dv{h}{x}(\hat x,u)\\
 \label{kalman:gain2}  L &=PC^TN^{-1} &\quad \dotex P &=AP+PA^T+M-PC^TN^{-1}CP.
 \end{align}
 The EKF has two main flaws when compared to the KF for time-invariant linear systems. First the linearized system around any trajectory is generally time-varying and the covariance matrix does not tend to a fixed value. Then, when $\hat x-x$ is large the linearized error equation can be a very erroneous approximation of the true error equation.

\section{Symmetry-preserving observers}
\label{sec:symmetry_preserving}
\subsection{Symmetry group of a system of differential equations}
Let $G$ be a group, and $M$ be a set. A group action can be defined on $M$ if to any $g\in G$ on can associate a diffeomorphic transformation $\phi_g:M\rightarrow M$ such that $\phi_{gh}=\phi_g\circ\phi_h$, and $(\phi_g)^{-1}=\phi_{g^{-1}}$, \emph{i.e.}, the group multiplication corresponds to the transformation composition, and the reciprocal elements correspond to reciprocal transformations.

\begin{defn}
G is a symmetry group of a system
of differential equations defined on $M$ if it maps
solutions to solutions. In this case we say the system is invariant.
\end{defn}

\begin{defn}
A vector field $w$ on $M$ is said
invariant if the system $\frac{d}{dt}z=w(z)$ is invariant.
\end{defn}

\begin{defn}
A scalar invariant  is a function $I:M\rightarrow \RR$ such that $I(\phi_g(z))=I(z)$ for all $g\in G$.
\end{defn}

Salar invariants and invariant vector fields can be built via Cartan's moving frame method \cite{olver-book99}.
\begin{defn}
A moving frame is a function $\gamma:M\rightarrow G$ such that $\gamma(\phi_g(z))=g\cdot\gamma(z)$ for all $g,z.$
\end{defn}

Suppose dim $G=r\leq $ dim $M$. Under some mild assumptions on the action (free, regular) there exists locally a moving frame. The sets $\mathcal O_z=\{\phi_g(z),~g\in G\}$ are called the group orbits. Let $K$ be a cross-section to the orbits. A moving frame  can be built locally via implicit functions theorem as the solution $g=\gamma(z)$ of the equation $\phi_g(z)=k$ where $k\in\mathcal O_z\cap K$. A complete set of functionnaly independent invariants is given by the non-constant  components of $\phi_{\gamma(z)}(z)$. Figure 1 illustrates those definitions and the moving frame method.

\begin{figure}
%\sidecaption
\includegraphics[scale=.27]{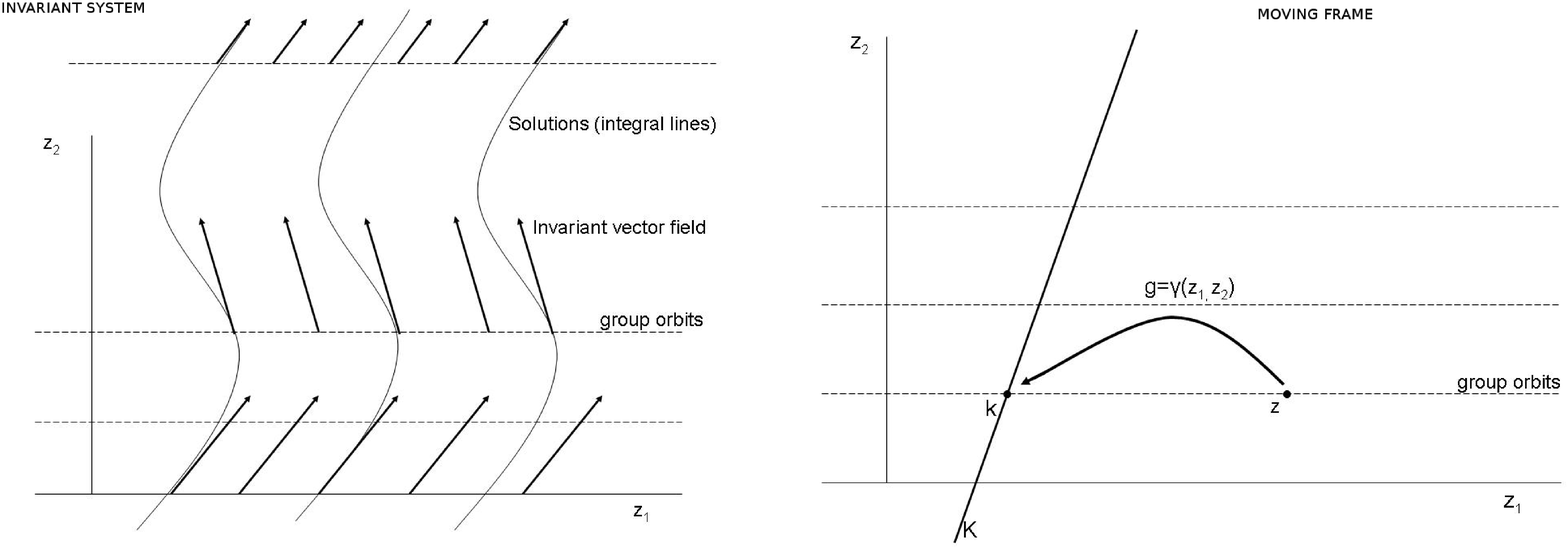}
\caption{\small{An illustrative example.   $M=\RR^2$, and the symmetry group is made of horizontal translations. We have $\phi_g(z_1,z_2)=(z_1+g,z_2)^T$ where $g\in G=\RR$. In local rectifying coordinates, every invariant system can be represented by a similar figure (under mild assumptions on the group action).
Left: Invariant system. The symmetry group maps each integral line of the vector field  into another integral line.  Right: Moving frame method.    $K$ is a cross-section to the orbits and $\gamma(z_1,z_2)$ is the group element that maps $(z_1,z_2)$ to $K$ along the orbit.  For exemple  if $K$ is the set $\{z_1\equiv 0\}$, the moving frame is $\gamma(z_1,z_2)=z_1$ and a complete set of invariants  is  $I(z_1,z_2)=z_2$.} }
\end{figure}

\subsection{Symmetry group of an observer}
Consider the general system \eqref{nonlin:eq1}.
Consider also the local group of transformations on $\XX\times\UU$
defined for any $x,u,g$ by \begin{align}\label{transfo:eq}
\phi_g(x,u)=\bigl(\varphi_g(x),\psi_g(u)\bigr),
\end{align}
where $\varphi_g$ and $\psi_g$ correspond to separate local group of transformations of $\XX$ and $\UU$.
\begin{prop} The system $\dotex x=f(x,u)$ is said invariant if it is invariant to the group action  \eqref{transfo:eq}.
\end{prop}
The group maps solutions to solutions if we have $\dotex X = f(X,U)$, where $(X,U)=(\varphi_g(x),\psi_g(u))$ for all $g\in G$. We understand from this definition, that $u$ can denote the control variables as usual, but it also denotes every feature of the environment that makes the system not behave the same way after it has been transformed (via $\varphi_g$). The action of $\psi_g$ is meant to allow some features of the environment to be also moved over. We would like the observer to be an invariant system for the  \emph{same} symmetry group.
\begin{defn}The observer \eqref{asymp:obs:eq}  is {invariant} or ``symmetry-preserving" if it is an invariant system for the group action $(\hat x,x, u, y)\mapsto \bigl(\varphi_g(x),\varphi_g(\hat x), \psi_g(u), h(\varphi_g( x), \psi_g(u))\bigr)$.
\end{defn}
In this case, the structure of the observer mimicks the nonlinear structure of the system. Let us recall how to build such observers (see \cite{arxiv-07} for more details). To do so, we need the output to be equivariant:
\begin{defn}The output is equivariant if  there exists  a group action on the output space  (via $\rho_g$) such that
$
h(\varphi_g(x),\psi_g(u))=\rho_g (h(x,u))
$ for all $g,x,u$.
\end{defn}
We will systematically assume the output is equivariant.  Let us define an invariant output error, instead of the usual linear output error $\hat
 y-y$:
\begin{defn} \label{inverr:defn}The smooth map $(\hat x,u,y)\mapsto E(\hat
x,u,y)\in\Rset^p$ is an \emph{invariant output error} if
\begin{itemize}
\item  $E\bigl(\varphi_g(\hat x),\psi_g(u),\rho_g(y)\bigr)
 =E(\hat x,u,y)$ for all~$\hat x,u,y$ (invariant)

\item  the map $y\mapsto E(\hat x,u,y)$ is invertible for all~$\hat
x,u$ (output)
\item  $E\bigl(\hat x,u,h(\hat x,u)\bigr)=0$ for all $\hat x,u$ (error)
\end{itemize}
\end{defn}
An invariant error is given (locally) by $E(\hat x,u,y)=\rho_{\gamma(\hat x,u)}(y)-\rho_{\gamma(\hat x,u)}(\hat y)$. Finally, an invariant frame $(w_1,...,w_n)$ on $\XX$, which  is a set of n linearly
point-wise independent invariant vector fields, i.e
$(w_1(x),...,w_n(x))$ is a basis of the tangent space to $\XX$ at
$x$. Once again such a frame can be built (locally) via the moving frame method.

\begin{prop}\cite{arxiv-07} The system $\dotex{\hat x}=F(\hat
x,u,y)$ is an invariant observer for the invariant
system $\dotex x=f(x,u)$ if and only if:
\begin{align}\label{inv:obs:eq}F(\hat x,u,y)=f(\hat x,u)
 +\sum_{i=1}^n\LL_i\bigl(I(\hat x,u),E(\hat x,u,y)\bigr)
 w_i(\hat x)\end{align}
where $E$ is an invariant output error, $I(\hat
x,u)$ is a complete set of scalar invariants, the
$\LL_i$'s are smooth functions such that for all~$\hat x$,
$\LL_i\bigl(I(\hat x,u),0\bigr)=0$, and $(w_1,...,w_n)$ is an
invariant frame.
\end{prop}
 The gains $\LL_i$ must be tuned in order to get some convergence properties if possible, and their magnitude should depend on the trade-off between measurement noise and convergence speed. The convergence analysis of the observer often relies on an invariant state-error:
\begin{defn} \label{invst:defn}The smooth map $(\hat x,x)\mapsto \eta(\hat
x,x)\in\Rset^n$ is an {invariant state error} if
$\eta(\varphi_g(\hat x),\varphi_g(x))=\eta(\hat x,x)$ (invariant),
the map $x\mapsto \eta(\hat x,x)$ is invertible for all~$\hat x$
(state), and $\eta(x,x)=0$ (error).
\end{defn}

\subsection{An example: symmetry-preserving observers for positive linear systems}
The linear system $\dotex x=Ax,~y=Cx$ admits scalings $G=\RR^*$  as a symmetry group via the group action $\phi_g(x)=g x$.  Every linear observer is obviously an invariant observer. The  unit sphere is a cross-section $K$ to the orbits. A moving frame maps the orbits to the sphere and thus writes $\gamma(x)=1/\norm x\in G$. A complete set of invariants is given locally by $n-1$ independent coordinates of $\phi_{\gamma(x)}(x)=x/\norm x$. Let $I(x)\in\RR^{n-1}$ be a complete set of independent invariants. $I(x)$ and $\norm x$ provide alternative coordinates named base and fiber coordinates. Moreover the system has a nice triangular structure in those coordinates. One can prove that $\dotex I(x(t))$ is an invariant function and thus it is necessarily of the form $g(I)$. As a result we have $\dotex I(x)=g(I(x))$ which does not depend on $\norm{x}$.

We have thus the following (general) result : if the restriction of the vector field on the cross-section is a contraction, it suffices to define a reduced observer on the orbits i.e. in our case a norm observer (which means that a scalar output suffices for observability). This is the case for instance when $A$ is a matrix whose coefficients are stricly positive (according to the Perron-Froebenius theorem). This fact was recently used in \cite{bonnabel-et-al:CDC11} to derive invariant asymptotic positive observers for positive linear systems.

\section{A first application to EKF SLAM}

\begin{figure}
%\sidecaption
\includegraphics[scale=.27]{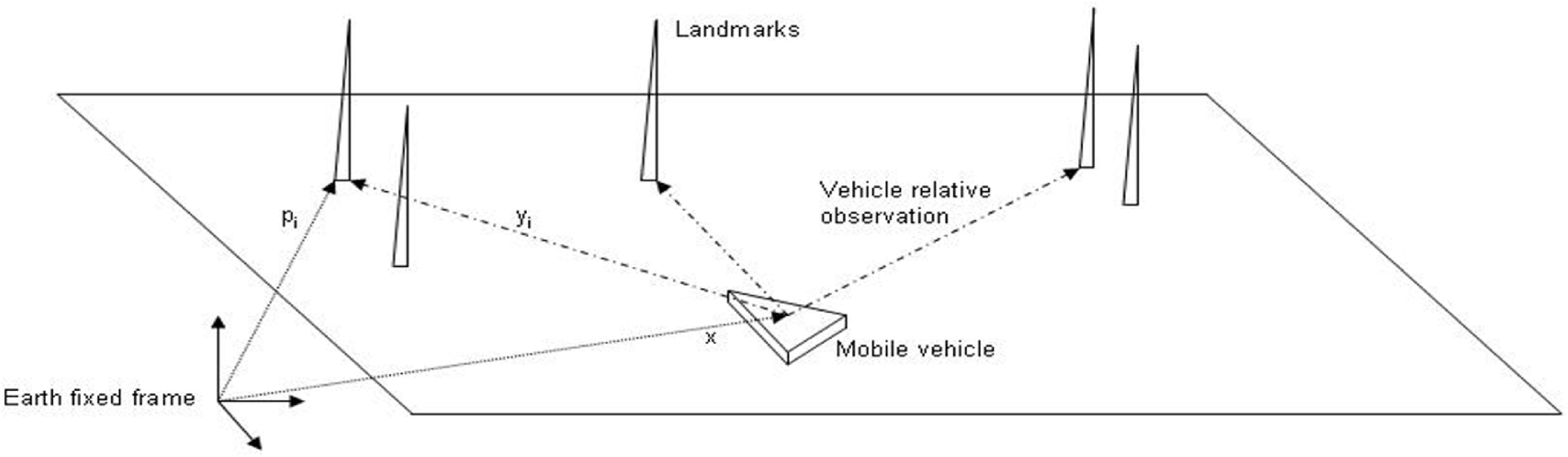}
\caption{\small{Vehicle taking relative measurements to environmental landmarks.} }
\end{figure}

Simultaneous localisation and mapping (SLAM) addresses the problem of building a map of an environment
from a sequence of sensor measurements obtained from a moving
robot. A solution to the SLAM problem has been seen for more than twenty years as a ``holy grail" in the robotics community since it would be a means to make a robot truly autonomous in an unknown environment.
A very well-known  approach that appeared in the early 2000's is the EKF SLAM  \cite{dissanayake-2001}.  Its main advantage is to formulate the problem in the form of a state-space model with additive
Gaussian noise and to provide convergence properties in the linear case (i.e. straight line motion). Indeed, the key idea is to include the position of the several landmarks (i.e. the map) in the state space.  This solution has been gradually replaced by other techniques such as FastSLAM, Graph SLAM etc.

In the framework of EKF SLAM, the problem of estimating online the trajectory of the robot as well as the location of all landmarks  without the need for any a priori
knowledge of location can be formulated as follows \cite{dissanayake-2001}. The vehicle state is defined by the position in the reference frame  (earth-fixed frame) $x\in\RR^2$  of
 the
centre of the rear axle and the  orientation of the vehicle axis $\theta$. The vehicle trusted motion relies on non-holonomic constraints. The landmarks are modeled as points and represented
by their position in the reference frame $p_i\in\RR^2$ where $1\leq i\leq N$. $u,v\in\RR$ are control inputs. Both
vehicle and landmark states are registered in the same frame of
reference. In a determistic setting (state noises turned off), the time evolution of the (huge) state vector is
\begin{align}\label{SLAM:eq}
\dot x&=u~R_\theta e_1,\quad
\dot\theta =uv,\quad
\dot p_i =0 \quad 1\leq i\leq N
\end{align}
where $e_1=(1,0)^T$ and $R_\theta$ is the rotation matrix of angle $\theta$.
Supposing that the data association between landmarks from one instant to the next is correctly done, the observation model  for the  $i$-th landmark  (disregarding measurement noise) is its position seen from the vehicle's frame:
$
z_i=R_{-\theta}(p_i-x)
$.
The standard EKF SLAM estimator has the form
\begin{align}\label{EKFSLAM:eq}
\dotex \hat x&=uR_{\hat \theta}e_1+L_x^k(\hat z_k- z_k), \quad
\dotex\hat \theta =uv+\sum_1^N L_\theta^k(\hat z_k- z_k),\quad
\dotex \hat p_i=\sum_1^N L_{i}^k(\hat z_k- z_k),\quad 1\leq i\leq N
\end{align}
where $\hat z_i=R_{-\hat\theta}(\hat p_i-\hat x)$ and where the $L_i$'s are the lines of $L$ tuned via the EKF equations \eqref{kalman:gain1}-\eqref{kalman:gain2}.

Here the group of symmetry of the system corresponds to Galilean invariances, and it  is made of rotations and translations of the plane $SE(2)$. Indeed, looking at Figure 2, it is obvious that the equations of motion are the same whether the first horizontal axis of the reference frame is pointing North, or East, or in any direction. For $g=(x_0,\theta_0)\in SE(2)$, the action of the group on the state space is $\varphi_g(x,\theta,p_i)=(R_{\theta_0}x+x_0,\theta+\theta_0,R_{\theta_0}p_i+x_0) $ and $\psi_g(u,v)=u,v$. The output is also unchanged by the group transformation as it is expressed in the vehicle frame and is thus insensitive to rotations and translations of the reference frame. Applying the theory of the last section, the observer above can be ``invariantized", yielding the following invariant observer:
\begin{align}\label{inv:ekf:eq}
\dotex \hat x&=uR_{\hat \theta}e_1+R_{\hat \theta}(\sum_1^N L_x^k(\hat z_k- z_k)), \quad
\dotex\hat \theta =uv+\sum_1^N L_\theta^k(\hat z_k- z_k),\quad
\dotex \hat p_i=R_{\hat \theta}(\sum_1^N L_{i}^k(\hat z_k- z_k))
\end{align}
It is easy to see that the invariant observer is much more meaningful, especially if the $L_i's$ are chosen as constant matrices \cite{arxiv-07}. Indeed, one could really wonder if it is sensible to correct vectors expressed in the reference frame directly  with measurements expressed in the vehicle frame. To be convinced, consider the following simple case : suppose $\hat\theta=\theta=\hat x=x=0$ remain fixed.   We have $\dotex (\hat p_i-p_i)=L_{i}(\hat p_i-p_i)$. Choosing $L_i=-k~I$ yields  $\dotex\norm{\hat p_i-p_i}^2=-k\norm{\hat p_i-p_i}^2$ leading to a correct estimation of landmark $p_i$. Now suppose that the vehicle has changed its orientation and $\hat \theta=\theta=\pi/2$.  The output error is now $R_{-\pi/2}(\hat p_i-p_i)$. With an observer of the form \eqref{EKFSLAM:eq} the same choice $L_i=-k~I$ yields $\dotex\norm{\hat p_i-p_i}=0$ and the landmark is not correctly estimated. On the other hand, with \eqref{inv:ekf:eq} we have in both cases $\dotex\norm{\hat p_i-p_i}^2=-k\norm{\hat p_i-p_i}^2$ ensuring convergence of $\hat p$ towards $p$.

Constant gains is a special (simple) choice, but the observer gains can also be tuned via Kalman equations. Indeed on can define  noises on the linearized invariant error system  and tune the $L_i$'s via Kalman equations (see Invariant EKF method  \cite{bonnabel-martin-salaun-cdc09}).
To sum up, any Luenberger  observer or EKF can be invariantized via equations \eqref{inv:ekf:eq}. This yields in the author's opinion  a much more meaningful non-linear observer that is well-adapted to the problem's structure. The invariantized observer \eqref{inv:ekf:eq} is simply a version of \eqref{EKFSLAM:eq} which is less sensitive to change of coordinates, and even if no proof can support this claim we believe it can only improve the performances of \eqref{EKFSLAM:eq}.

\section{Particular case where the state space coincides with its symmetry group}

Over the last half decade, invariant observers on Lie groups for low-cost aided inertial navigation have been studied by several teams in the world,  \cite{mahony-et-al-IEEE,arxiv-07,vasconcelos} to name a few. Several powerful  convergence results have been obtained. They are all linked to the special properties of the invariant state error on a Lie group. To recap briefly the construction of invariant observers on Lie groups \cite{arxiv-08}, we assume that the symmetry group $G$ is a matrix group, and that $\XX=G$. The system is assumed to be invariant to left multiplications i.e. $
\dotex X=X\Omega(t).$  We have indeed for any $g\in G$ that $\dotex (gX)=(gX)\Omega$.
For instance the motion of the vehicle in the considered SLAM problem   $\dot x=uR_\theta e_1,\dot\theta=uv$ can be viewed as a left-invariant system on the Lie group SE(2) via the matrix representation $$
X=\begin{pmatrix}R_\theta & x\\0_{1\times 2}&1\end{pmatrix},\quad \Omega=\begin{pmatrix}\omega_x & ue_1\\0_{1\times 2}&0\end{pmatrix},~\text{with}~\omega_x=\begin{pmatrix}0 & -uv\\uv&0\end{pmatrix}$$
Suppose the output $y=h(X)$ is equivariant, i.e. there exists a group action on the output space such that $h(gX)=\rho_g(X)$. In this case the invariant observer  \eqref{inv:obs:eq} can be written intrinsically
$$
\dotex \hat X=\hat X\Omega+\hat X L(\rho_{\hat X^{-1}}(y)).
$$with $L(e)=0$ where $e$ is the group identity element. The invariant state error is the natural group difference $\eta=X^{-1}\hat X$ and the error equation is $$\dotex \eta=[\Omega, \eta]+\eta L\circ h(\eta^{-1})$$ A remarkable fact is that the error equation only depends on $\eta$ and $\Omega$, whereas the system is non-linear and the error should also depend on $\hat X$ (think about the EKF which is based on a linearization around any $\hat X$ at each time). Moreover, if $\Omega=cst$, the error equation is clearly \emph{autonomous}. Thus the motion primitives generated by constant $\Omega$ are special trajectories called ``permanent trajectories". Around such trajectories one can always achieve local convergence (as soon as the linearized system is observable).

It is worth noting this property was recently used to derive a non-linear \emph{separation principle} on Lie groups \cite{bonnabel-et-al:ifac11}. It applies to some cart-like underactuated vehicles and some underwater or aerial fully actuated vehicles.

An even more interesting case occurs when the output satisfies right-equivariance i.e.,  $h(Xg)=\rho_g(h(X))$. In this case we let the input be $u=\Omega$ and we consider the  action of $G$ by right multiplication, i.e. $\varphi_g(X)=Xg$ and $\psi_g(\Omega)= g^{-1}\Omega g$. The output is equivariant as $h(\varphi_g(X))=\rho_g\circ h(X)$. The invariant observer associated with this group of symmetry writes
$
\dotex \hat X=\hat X\Omega+ L(\rho_{\hat X^{-1}}(y))\hat X
$. The invariant state error is $\eta=\hat X X^{-1}$ and the error equation is
\begin{align}\label{rerror:eq}\dotex\eta=\hat X\Omega X^{-1}+L(h(\eta^{-1}))\eta-\hat X\Omega X^{-1}=L(h(\eta^{-1}))\eta\end{align}The error equation is completely autonomous ! In particular the linearized system around \emph{any} trajectory is the same time-invariant system. Autonomy is the key for numerous powerful convergence results for observers on Lie groups  see e.g. \cite{lageman10,vasconcelos,arxiv-08}.

\section{A new result in EKF SLAM}

In this section we propose a new non-linear observer for EKF SLAM with guaranteed convergence properties.
In the SLAM problem  the state space is much bigger than its symmetry group. The orbits have dimension 3 and thus there are $N+1-3+2=N$ invariants (dimension of the cross-section, see Fig.1). Thus an autonomous error equation seems to be out of reach. Suprinsingly considering the symmetry group of rotations and translations in the vehicle frame yields such a result. A simple trick makes it obvious. Consider the following matrix representation:
$$
X=\begin{pmatrix}R_\theta & x\\0_{1\times 2}&1\end{pmatrix},\quad P_i=\begin{pmatrix}R_\theta & p_i\\0_{1\times 2}&1\end{pmatrix}, \quad \Omega=\begin{pmatrix}\omega_x & ue_1\\0_{1\times 2}&0\end{pmatrix},\quad \Omega_i=\begin{pmatrix}\omega_x & 0\\0_{1\times 2}&0\end{pmatrix}$$
The equations of the system \eqref{SLAM:eq} can be written
$
\dotex X=X\Omega,~ \dotex P_i=P_i\Omega_i,~1\leq i\leq N
$
and the system can be viewed as a left-invariant dynamics system on the (huge) Lie group $G\times\cdots\times G$. Let $\eta_x=\hat XX^{-1},\eta_i=\hat P_iP_i^{-1}$ be the invariant state error. The system has the invariant output errors
$
\tilde Y_i=R_{\hat\theta}(\hat z_i-z_i)$, i.e. $\begin{pmatrix}\tilde Y_i&1\end{pmatrix}^T=(\eta_i-\eta_x)H
$ for $1\leq i\leq N$ where  $H=\begin{pmatrix}0_{1\times 2}&1\end{pmatrix}^T$.
 Consider the following invariant observer
$
\dotex \hat X=\hat X\Omega+L_X(\tilde Y_1,\cdots, \tilde Y_N)\hat X,\quad \dotex \hat P_i=\hat P_i\Omega_i+L_i(\tilde Y_1,\cdots, \tilde Y_N)\hat P_i
$.
From \eqref{rerror:eq}, the (non-linear) error equation is completely autonomous reminding the linear case \eqref{err:eq3}. It implies the following global convergence result for the non-linear deterministic system:

\begin{prop}
Consider the SLAM problem \eqref{SLAM:eq} without noise. The following observer
\begin{align*}
\dotex\hat \theta =uv,\quad \dotex \hat x&=uR_{\hat \theta}e_1, \quad
\dotex \hat p_i=k_i~R_{\hat\theta}(\hat z_i-z_i)
\end{align*}
with $k_i>0$ is such that $\dotex(R_{\hat\theta}(\hat z_i-z_i))=-k_i~R_{\hat\theta}(\hat z_i-z_i)$, i.e., all the estimation errors
 $
(\hat z_i-z_i),~ 1\leq i\leq N$
converge globally exponentially to zero with rate $k_i$, which means  the vehicle trajectory and the map are correctly identified. The parameter $k_i$ must be tuned according to the level of noise associated to landmark $i$, and vehicle sensors' noise.
\end{prop}
If one wants to define noise covariance matrices $M,N$ to tune the observer (and compute an estimation $P$  of the covariance error matrix at each time), it is also possible to define a modified EKF with guaranteed convergence properties:
\begin{prop}
Consider the SLAM problem \eqref{SLAM:eq}. Let $E=(R_{\hat\theta}(\hat z_i-z_i))_{1\leq i\leq N}$ be the invariant output error.  Let $e_3$ be the vertical axis. Consider the observer
\begin{align*}\label{inv:ekf2:eq}
\dotex\hat \theta =uv+\mathcal L_\theta(E),\quad \dotex \hat x&=uR_{\hat \theta}e_1+\mathcal L_\theta(E) e_3\wedge \hat x+ \mathcal L_x(E), \quad
\dotex \hat p_i=\mathcal L_\theta(E) e_3\wedge \hat p+ \mathcal L_i(E)
\end{align*}
Let $\eta=(\tilde\theta,\tilde x,\tilde p_1,\cdots,\tilde p_n)$ be the invariant state error where $
\tilde\theta=\hat\theta-\theta,~\tilde x=\hat x-R_{\tilde\theta} x,~\tilde p_i=\hat p_i-R_{\tilde\theta} p_i$. The state error equation is autonomous, i.e. $\dotex\eta$  only depends on $\eta$. It is thus completely independent of the trajectory and of $u(t),v(t)$. The linearized error equation writes
 $
\dotex \delta\eta=(LC)\delta \eta
$ where $L$ can be freely chosen and  $C$ is a fixed matrix.
As in the usual EKF method, one can define covariance matries $M$, $N$, build a Kalman filter for the linearized system, i.e.  tune $L$ via the usual equations \eqref{kalman:gain2} i.e. $\dot P=M-PC^TN^{-1}CP$,  $L=PC^TN^{-1}$, and implement it on the non-linear model. All the convergence results on $P$ and $L$ valid for stationnary systems \eqref{eq:lin_sys} with $A=0,~B=0, ~D=0$ apply.
\end{prop}

\begin{figure}
%\sidecaption
\includegraphics[scale=.42]{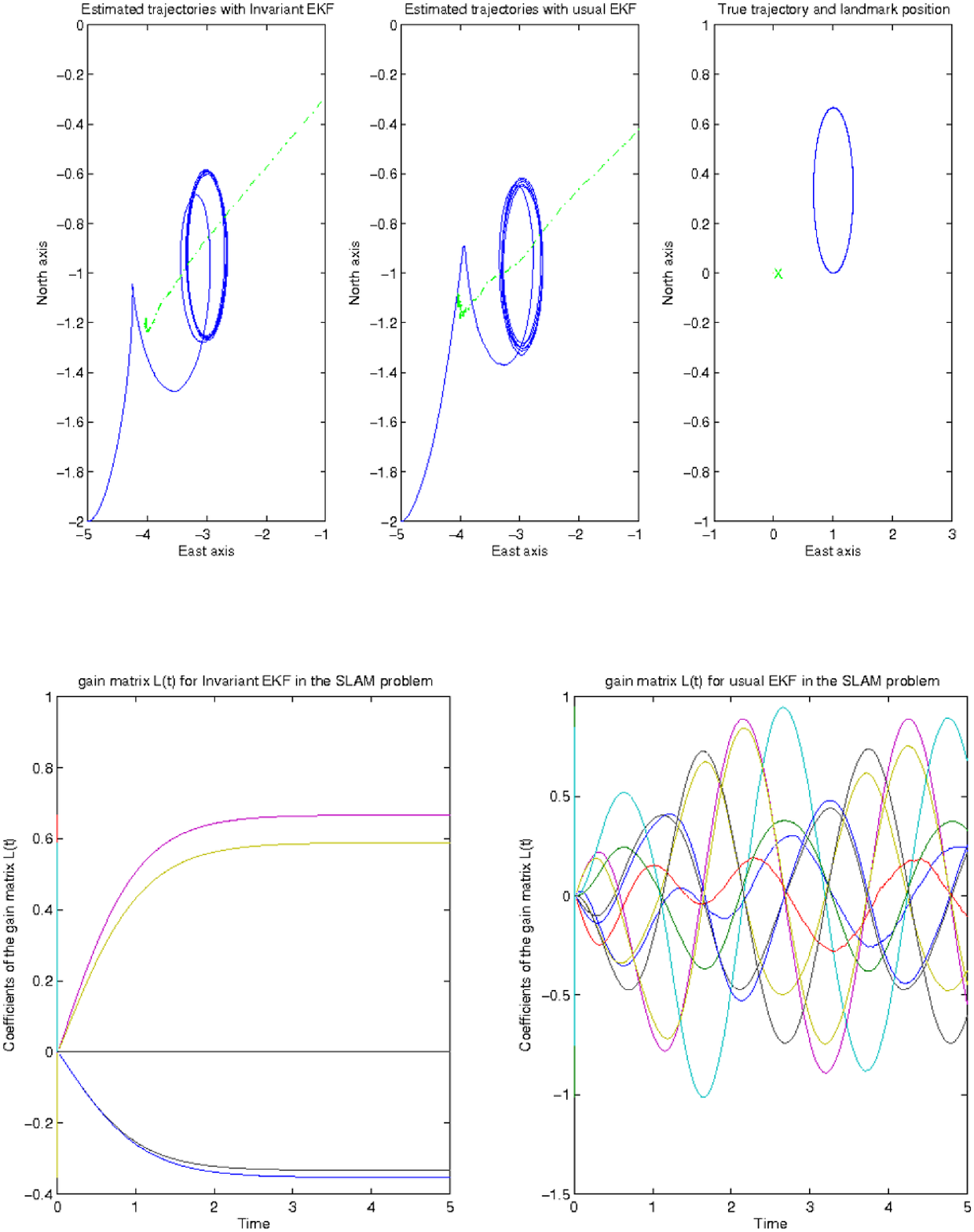}
\caption{\small{Simulations with one landmark and a car moving over a
circular path with a $20\%$ noise. Up: 1-Estimated vehicle
trajectory (plain blue line) and landmark position (dashed green
line) with Invariant EKF, 2-Estimation with the usual EKF, 3- true
vehicle trajectory (plain blue line) and landmark position (green
cross). After a short transient, the trajectory is correctly
identified for both observers (up to a rotation-translation). Bottom
: 1-coefficients of $L(t)$ over time for Invariant EKF,
2-coefficients of $L(t)$ for EKF.  Wee see the EKF gain matrix is
permanently adapting to the motion of the car (right) whereas its
invariant counterpart (left) is directly expressed in well-adapted
variables.}}\label{simuf}
\end{figure}

Simulations (Fig. \ref{simuf}) with one landmark and noisy
measurements indicate the modified EKF (IEKF) behaves very similarly,
or slightly better than the EKF, but the gain matrix tends quickly
to a fixed matrix $L$  independently from the trajectory and the
inputs $u,v$. So the Invariant EKF proposed in this paper 1- is
incomparably cheaper computationaly as it relies on a constant
matrix $L$ that can be computed offline once and for all (the number
of landmarks can thus be much increased) 2- is such that the
linearized error system is stable as soon as $LC$ has negative
eigenvalues, which is easy to verify.

\begin{rem}The calculations above are valid on $SE(3)$ and the results apply to 6 DOF SLAM. 
\end{rem}

\begin{thebibliography}{10}

\bibitem{aghannan-rouchon-cdc02}
N.~Aghannan and P.~Rouchon.
\newblock On invariant asymptotic observers.
\newblock In {\em 41st IEEE Conference on Decision and Control}, pages
  1479--1484, 2002.

\bibitem{bonnabel-martin-salaun-cdc09}
S.~Bonnabel, P.~Martin, and E.~Sala{u}n.
\newblock Invariant extended kalman filter: Theory and application to a
  velocity-aided attitude estimation problem.
\newblock In {\em IEEE Conference on Decision and Control}, 2009.

\bibitem{arxiv-07}
S.~Bonnabel, Ph. Martin, and P.~Rouchon.
\newblock Symmetry-preserving observers.
\newblock {\em IEEE Trans. on Automatic Control}, 53(11):2514--2526, 2008.

\bibitem{arxiv-08}
S.~Bonnabel, Ph. Martin, and P.~Rouchon.
\newblock Non-linear symmetry-preserving observers on lie groups.
\newblock {\em IEEE Trans. on Automatic Control}, 54(7):1709 -- 1713, 2009.

\bibitem{bonnabel-et-al:ifac11}
S.~Bonnabel, Ph. Martin, P.~Rouchon, and E.~Salaun.
\newblock A separation principle on lie groups.
\newblock In {\em IFAC (available on Arxiv)}, 2011.

\bibitem{bonnabel-et-al:CDC11}
S.~Bonnabel and R.~Sepulchre.
\newblock Contraction and observer design on cones.
\newblock {\em Arxiv}, 2011.

\bibitem{bullo-murray-auto99}
F.~Bullo and R.M. Murray.
\newblock Tracking for fully actuated mechanical systems: A geometric
  framework.
\newblock {\em Automatica}, 35(1):17--34, 1999.

\bibitem{dissanayake-2001}
G.~Dissanayake, P.~Newman, H.F. Durrant-Whyte, S.~Clark, and M.~Csobra.
\newblock A solution to the simultaneous localisation and mapping (slam)
  problem.
\newblock {\em IEEE Trans. Robot. Automat.}, 17:229--241, 2001.

\bibitem{grizzle-marcus-ieee85}
J.W. Grizzle and S.I. Marcus.
\newblock The structure of nonlinear systems possessing symmetries.
\newblock {\em IEEE Trans. Automat. Control}, 30:248--258, 1985.

\bibitem{lageman10}
C.~Lagemann, J.~Trumpf, and R.~Mahony.
\newblock Gradient-like observers for invariant dynamics on a lie group.
\newblock {\em IEEE Trans. on Automatic Control}, 55:2:367 -- 377, 2010.

\bibitem{mahony-et-al-IEEE}
R.~Mahony, T.~Hamel, and J-M Pflimlin.
\newblock Nonlinear complementary filters on the special orthogonal group.
\newblock {\em IEEE-Trans. on Automatic Control}, 53(5):1203--1218, 2008.

\bibitem{martin-et-al-cocv03}
Ph. Martin, P.~Rouchon, and J.~Rudolph.
\newblock Invariant tracking.
\newblock {\em ESAIM: Control, Optimisation and Calculus of Variations},
  10:1--13, 2004.

\bibitem{morin-ieee03}
P.~Morin and C.~Samson.
\newblock Practical stabilization of driftless systems on lie groups, the
  transverse function approach.
\newblock {\em IEEE Trans. Automat. Control}, 48:1493--1508, 2003.

\bibitem{olver-book99}
P.~J. Olver.
\newblock {\em Classical Invariant Theory}.
\newblock Cambridge University Press, 1999.

\bibitem{vasconcelos}
J.F. Vasconcelos, R.~Cunha, C.~Silvestre, and P.~Oliveira.
\newblock A nonlinear position and attitude observer on se(3) using landmark
  measurements.
\newblock {\em Systems Control Letters}, 59:155--166, 2010.

\end{thebibliography}
\end{document}